\documentclass[11pt]{amsart}
\usepackage{amsmath,amssymb,latexsym,amsfonts,amscd}
\newtheorem{theorem}{Theorem}[section]
\newtheorem{proposition}[theorem]{Proposition}
\newtheorem{definition}[theorem]{Definition}

\newtheorem{lemma}[theorem]{Lemma}
\newtheorem{corollary}[theorem]{Corollary}

\begin{document}

\title[Switala's Matlis duality]{On Switala's Matlis duality for $\mathcal D$-modules}
\author{Gennady Lyubeznik}
\address{Dept. of Mathematics, University of Minnesota, Minneapolis,
MN 55455}
\email{gennady@math.umn.edu}
\begin{abstract}We show that N. Switala's Matlis duality for $\mathcal D$-modules yields a generalization of a result of R. Hartshorne and C. Polini.
\end{abstract}
\keywords{complete local ring, holonomic module, de Rham cohomology}
\subjclass[2000]{Primary 13H05, 13N10, 16S32}
\thanks{NSF for support through grant DMS-1500264 is gratefully acknowledged}

\maketitle

\section{Introduction}

Let $(R,\mathfrak m)$ be a commutative Noetherian complete equicharacteristic local ring and let $\mathcal D=\mathcal D(R,k)$ be the ring of $k$-linear differential operators of $R$, where $k\subset R$ is a fixed coefficient field of $R$ ($k$ exists as $R$ is complete). 

Let $E$ be the injective hull of $R/\mathfrak m\cong k$ in the category of $R$-modules. Matlis duality $D(-)={\rm Hom}_R(-,E)$ is a contravariant exact functor from $R$-modules to $R$-modules. N. Switala \cite[Section 4]{S1} has shown that if $M$ is a left (resp. right) $\mathcal D$-module, then $D(M)$ carries the structure of a right (resp. left) $\mathcal D$-module. In this short note our main result is the following.

\begin{theorem}\label{main}
Let $M$ be a left $\mathcal D$-module and let $N$ be a right $\mathcal D$-module. There is a $k$-linear bijection $${\rm Ext}^i_{\mathcal D{\rm  -mod}}(M,D(N))\cong {\rm Ext}^i_{{\rm mod-}\mathcal D}(N,D(M))$$
which is functorial in $M$ and $N$, where $\mathcal D{\rm -mod}$ and ${\rm mod-}\mathcal D$ denote the categories of, respectively, left and right $\mathcal D$-modules.
\end{theorem}

R. Hartshorne and C. Polini proved the following important result \cite[5.2]{HP}.

\medskip

\begin{theorem}\label{HP}
Notation being as above, assume the characteristic of $k$ is 0 and $R=k[[x_1,\dots,x_n]]$ is the ring of formal power series in n variables over $k$. Let $M$ be a left $\mathcal D$-module and let $H^i_{\rm dR}(M)$ denote the $i$-th de Rham cohomology group of $M$. If $M$ is holonomic, then $${\rm dim}_kH^n_{\rm dR}(M)={\rm dim}_k{\rm Hom}_{\mathcal D{\rm -mod}}(M,E).$$
\end{theorem}

As an application of our Theorem \ref{main} we show that Switala's Matlis duality theory yields the following generalization of Theorem \ref{HP}.

\begin{theorem}\label{maincoro}
Notation being as in Theorem \ref{HP}, for every $i$ $${\rm dim}_kH^{n-i}_{\rm dR}(M)={\rm dim}_k{\rm Ext}^i_{\mathcal D{\rm -mod}}(M,E)$$
\end{theorem}

Theorem \ref{HP} is the $i=0$ case of Theorem \ref{maincoro}.

\section{A review of Switala's Matlis duality}

Let $L$ be an $R$-module. The $k$-module Hom$_k(L,k)$ has a structure of $R$-module as follows: if $f:L\to k$ is a $k$-linear map and $r\in R$, the map $g=rf:L\to k$ is defined by $g(x)=f(rx)$ for every $x\in L$. The following definition is equivalent to but somewhat simpler than Switala's original definition \cite[Definitions 3.13, 3.14(c)]{S1}.

\begin{definition}
A $k$-linear map $f:L\to k$ is {\it $\Sigma$-continuous} if for every $x\in L$ there is some $i\geq 0$ so that $f(\mathfrak m^ix)=0$, where $\mathfrak m\subset R$ is the maximal ideal. The set of all $\Sigma$-continuous maps $L\to k$ is denoted by $D^\Sigma(L)$.
\end{definition}

If $f:L\to k$ is $\Sigma$-continuous, then so is $rf:L\to k$ because $r\mathfrak m^i\subset \mathfrak m^i$ and therefore $f(\mathfrak m^ix)=0$ implies $f(r\mathfrak m^ix)=0$. Thus $D^\Sigma(L)$ is an $R$-submodule of Hom$_k(L,k)$.

If $\phi:L'\to L$ is an $R$-module homomorphism and $f:L\to k$ is $\Sigma$-continuous, then the composition $f\circ\phi:L'\to k$ also is $\Sigma$-continuous. This defines a map $D^\Sigma(\phi):D^\Sigma(L)\to D^\Sigma(L')$ that sends $f\in D^\Sigma(L)$ to $f\circ\phi\in D^\Sigma(L')$. It is not hard to see that $D^\Sigma(\phi)$ is an $R$-module homomorphism. 

Thus one gets a contravariant functor $$D^\Sigma:R{\rm -mod}\to R{\rm -mod}$$ which takes every $R$-module $L$ to $D^\Sigma(L)$ and every $R$-module homomorphism $\phi: L\to L'$ to $D^\Sigma(\phi):D^\Sigma(L')\to D^\Sigma(L).$

Switala has shown \cite[Theorem 3.15]{S1} that there is an isomorphism of functors $$D^\Sigma\cong D$$ where $D={\rm Hom}_R(-, E)$ is the classical Matlis duality functor. Thus as long as one sticks to $R$-modules and $R$-module homomorphisms, $D^\Sigma$ does not give anything new compared to $D$. 

But $D^\Sigma$ makes it possible to dualize more maps between $R$-modules, than just $R$-module homomorphisms. For example, assume $L$ is a (left or right) $\mathcal D$-module, $\delta\in \mathcal D$ is a differential operator and $\delta:L\to L$ is the action of $\delta$ on $L$. Switala shows that if $f:L\to k$ is $\Sigma$-continuous, then the composition $f\circ\delta:L\to k$ also is $\Sigma$-continuous. This makes it possible to define the action of $\delta$ on $D^\Sigma(L)$, namely, $\delta(f)=f\circ\phi$ for every $f\in D^\Sigma$ \cite[Corollary 4.10]{S1}. 

This definition of the action of differential operators on $D^\Sigma$ gives $D^\Sigma$ a structure of $\mathcal D$-module. More precisely, if $L$ is a left $\mathcal D$-module, then $D^\Sigma(L)$ becomes a right $\mathcal D$-module and if $L$ is a right $\mathcal D$-module, then $D^\Sigma(L)$ becomes a left $\mathcal D$-module and if $\phi:L\to L'$ is a map of left (resp. right) $\mathcal D$-modules, then the map $D^\Sigma(\phi):D^\Sigma(L')\to D^\Sigma(L)$ is a map of right (resp. left) $\mathcal D$-modules. Thus $D^\Sigma$ induces two contravariant exact functors on $\mathcal D$-modules:
$$D^\Sigma:\mathcal D{\rm -mod}\to {\rm mod-}\mathcal D$$
$$D^\Sigma:{\rm mod-}\mathcal D\to \mathcal D{\rm -mod},$$
where $\mathcal D{\rm -mod}$ and ${\rm mod-}\mathcal D$ are the categories of left and right $\mathcal D$-modules respectively \cite[Proposition 4.11]{S1}.

\section{Proof of Theorem \ref{main}}

Given the $R$-module structures on Hom$_k(M,k)$ and Hom$_k(N,k)$ defined in the beginning of the preceding section, there are $R$-linear maps $$\tau(M,N):{\rm Hom}_R(M,{\rm Hom}_k(N,k))\to {\rm Hom}_R(N,{\rm Hom}_k(M,k))$$  $$\tau(N,M):{\rm Hom}_R(N,{\rm Hom}_k(M,k))\to {\rm Hom}_R(M,{\rm Hom}_k(N,k))$$ defined as follows. If $f\in {\rm Hom}_R(M,{\rm Hom}_k(N,k))$, then $g=\tau(M,N)(f)\in {\rm Hom}_R(N,{\rm Hom}_k(M,k))$ is defined by the formula
\begin{equation}\label{E:first}
g(y)(x)=f(x)(y)
\end{equation}
for every $x\in M$ and $y\in N$, and, similarly, if $g\in {\rm Hom}_k(N,{\rm Hom}_k(M,k))$, then $f=\tau(N,M)(g)\in {\rm Hom}_k(M,{\rm Hom}_k(N,k))$ is defined by the same formula (\ref{E:first}). Clearly, $g=\tau(M,N)(f)$ if and only if $f=\tau(N,M)(g)$, i.e. the two maps $\tau(M,N)$ and $\tau(N,M)$ are inverses of each other and therefore establish an $R$-linear bijection $${\rm Hom}_R(M,{\rm Hom}_k(N,k))\cong {\rm Hom}_R(N,{\rm Hom}_k(M,k)).$$

Let $f\in {\rm Hom}_R(M,{\rm Hom}_k(N,k))$ and $g\in {\rm Hom}_R(N,{\rm Hom}_k(M,k))$ correspond to each under this bijection, i.e. $g=\tau(M,N)(f)$ and $f=\tau(N,M)(g)$. We claim that  $f(x):N\to k$ is $\Sigma$-continuous for every $x\in M$ if and only if $g(y):M\to k$ is $\Sigma$-continuous for every $y\in N$. Indeed, assume $f(x)$ is $\Sigma$-continuous for every $x\in M$. Pick any $x\in M$ and $y\in N$. Since $f(x)\in D^\Sigma(N)$, there exists $i$ such that $f(x)(\mathfrak m^iy)=0$. Let $r\in \mathfrak m^i$. By definition, $(rf(x))(y)=f(x)(ry)$, hence $(rf(x))(y)=0$. Since $f$ is a map of $R$-modules, $rf(x)=f(rx)$. Hence $f(rx)(y)=0$, i.e. $g(y)(rx)=0$. Since this holds for every $r\in \mathfrak m^i$, the map $g(y)$ is $\Sigma$-continuous for every $y$. A similar argument shows that if $g(y)$ is $\Sigma$-continuous for every $y\in N$, then $f(x)$ is $\Sigma$-continuous for every $x\in M$. This proves the claim. 

Thus the maps $\tau(M,N)$ and $\tau(N,M)$ (once they are restricted to \hbox{the $R$-sub-} modules ${\rm Hom}_R(M,D^\Sigma(N))$ of ${\rm Hom}_R(M,{\rm Hom}_k(N,k))$ and ${\rm Hom}_R(N,D^\Sigma(M))$ of ${\rm Hom}_R(N,{\rm Hom}_k(M,k))$) establish an $R$-linear bijection $${\rm Hom}_R(M,D^\Sigma(N))\cong {\rm Hom}_R(N,D^\Sigma(M)).$$

Let $f\in {\rm Hom}_R(M,D^\Sigma(N))$ and $g\in {\rm Hom}_R(N,D^\Sigma(M))$ correspond to each under this bijection, i.e. $g=\tau(M,N)(f)$ and $f=\tau(N,M)(g)$. We claim that $f$ is a map of left $D$-modules if and only if $g$ is a map of right $D$-modules. Indeed, $f$ is a map of left $\mathcal D$-modules if and only if $f(\delta x)=\delta(f(x))$ for every $x\in M$ and every $\delta\in \mathcal D$. Set $h=f(x):N\to k$. By definition $\delta(h):N\to k$ is defined by $(\delta(h))(y)=h(y\delta)$ for every $y\in N$. Hence $f$ is a map of left $\mathcal D$-modules if and only if $$f(\delta x)(y)=f(x)(y\delta)$$ for every $x\in M, y\in N$ and $\delta\in \mathcal D$. A similar argument shows that $g$ is a map of right $\mathcal D$-modules if and only if $$g(y\delta)(x)=g(y)(\delta x)$$ for all $x\in M, y\in N$ and $\delta\in \mathcal D$. Since by formula (\ref{E:first}), $f(\delta x)(y)=g(y)(\delta x)$ and $f(x)(y\delta)=g(y\delta)(x)$, the claim is proven. 

Thus the maps  $\tau(M,N)$ and $\tau(N,M)$ (after restriction to the $k$-submodules ${\rm Hom}_{\mathcal D{\rm -mod}}(M,D^\Sigma(N))$ of ${\rm Hom}_R(M,D^\Sigma(N))$ and ${\rm Hom}_{{\rm mod-}\mathcal D}(N,D^\Sigma(M))$ of ${\rm Hom}_R(N,D^\Sigma(M))$) establish a $k$-linear bijection $${\rm Hom}_{\mathcal D{\rm -mod}}(M,D^\Sigma(N))\cong {\rm Hom}_{{\rm mod-}\mathcal D}(N,D^\Sigma(M)).$$
Since this bijection is clearly functorial in $M$ and $N$, the $i=0$ case of the theorem is proven.

We need the following lemma.

\begin{lemma}\label{lemma}
If $P$ is a projective left (resp. right) $\mathcal D$-module, then $D(P)$ is an injective right (resp. left) $\mathcal D$-modules. 
\end{lemma}

\emph{Proof of the lemma}. Assume $P$ is projective left $\mathcal D$-module. Let $$0\to N'\to N\to N''\to 0$$ be an exact sequence of right $\mathcal D$-modules. The induced sequence $$0\to {\rm Hom}_{{\rm mod-}\mathcal D}(N'', D(P))\to {\rm Hom}_{{\rm mod-}\mathcal D}(N, D(P))\to {\rm Hom}_{{\rm mod-}\mathcal D}(N', D(P))\to 0$$ is exact because by the $i=0$ case of the theorem that has just been proven, this sequence is isomorphic to the sequence $$0\to {\rm Hom}_{\mathcal D{-mod}}(P, D(N''))\to {\rm Hom}_{\mathcal D{-mod}}(P, D(N))\to {\rm Hom}_{\mathcal D{-mod}}(P, D(N'))\to 0$$ which is exact since $P$ is projective and $D$ is an exact functor. Therefore $D(P)$ is injective. 

For a projective right $\mathcal D$-module $P$ the proof is the same except "left" must be replaced by "right" and vice versa.\qed

\medskip

We continue with the proof of Theorem \ref{main}. Let $$\dots\to P_1\to P_0\to M\to 0$$ be a projective resolution of $M$ in the category of left $\mathcal D$-modules. The $k$-module Ext$^i_{\mathcal D{\rm -mod}}(M, D(N))$ is the $i$th cohomology of the induced complex  $$0\to {\rm Hom}_{\mathcal D{\rm -mod}}(P_0, D(N))\to {\rm Hom}_{\mathcal D{\rm -mod}}(P_1, D(N))\to \dots$$ By the $i=0$ case of the theorem that has already been proven, this complex is isomorphic to the complex $$0\to {\rm Hom}_{{\rm mod-}\mathcal D}(N, D(P_0))\to {\rm Hom}_{{\rm mod-}\mathcal D}(N, D(P_1))\to \dots$$ The $k$-module Ext$^i(N, D(M))$ is the $i$-th cohomology of this complex since by Lemma \ref{lemma} $D(P^\bullet)$ is an injective resolution of $D(M)$. This completes the proof of Theorem \ref{main}. \qed

\medskip

The natural action of $\mathcal D$ on $R$ makes $R$ a left $\mathcal D$-module. Hence by Switala's Matlis duality, $E=D(R)$ is a right $\mathcal D$-module. 

\begin{corollary}\label{RE}
If $N$ is a right $\mathcal D$-module, there exist $k$-linear bijections $${\rm Ext}^i_{\mathcal D{\rm  -mod}}(R,D(N))\cong {\rm Ext}^i_{{\rm mod-}\mathcal D}(N,E)$$
that are functorial in $N$.
\end{corollary}

\emph{Proof.} This is immediate from Theorem \ref{main} upon setting $M=R$.\qed

\section{Proof of Theorem \ref{maincoro}}

Throughout this section the field $k$ has characteristic 0 and $R$ is regular, i.e. $R=k[[x_1,\dots,x_n]]$ is the ring of formal power series in $n$ variables over $k$. Let $$\partial _i=\frac{\partial}{\partial x_i}:R\to R$$ be the $k$-linear partial differentiation with respect to $x_i$. The ring $\mathcal D$ is the free left (as well as right) $R$-module on the monomials $\partial_1^{a_1}\cdots\partial_n^{a_n}$. There is a ring anti-automorphism $$t:\mathcal D\to \mathcal D$$ (that is an isomorphism of the underlying additive groups that reverses the order of multiplication, i.e.  $t(\delta_1\delta_2)=t(\delta_2)t(\delta_1)$ for all $\delta_1,\delta_2\in \mathcal D$) defined by $$t(r\partial_1^{a_1}\cdots\partial_n^{a_n})=(-\partial_1)^{a_1}\cdots(-\partial)_n^{a_n}r$$ for all $r\in R$. This anti-automorphism $t$ is called transposition.

We can view a right $\mathcal D$-module $M$ as a left $\mathcal D$-module $t(M)$ via this transposition anti-automorphism \cite[Definition 4.14]{S1}, i.e. $$\delta*x=xt(\delta)$$ for every $x\in M$ and $\delta\in \mathcal D$ where $*$ denotes the left $\mathcal D$-action on $t(M)$. Thus $t(M)$ has the same elements and the same $R$-module structure as $M$ but the $D$-module structure is different. If $\phi:M\to M'$ is a map of right $\mathcal D$-modules, the map $t(\phi):t(M)\to t(M')$, which is the same map as $\phi$ on the underlying sets, is a map of left $\mathcal D$-modules. 

The resulting functor from right $\mathcal D$-modules to left $\mathcal D$-modules that sends a right $\mathcal D$-module $M$ to the left $\mathcal D$-module $t(M)$ and sends a map of right $\mathcal D$-modules $\phi:M\to M'$ to the map of left $\mathcal D$-modules $t(\phi):t(M)\to t(M')$ is covariant, exact, induces bijections on Hom-sets (i.e . ${\rm Hom}_{{\rm mod-}D}(M,M')\cong{\rm Hom}_{\mathcal D{\rm -mod}}(t(M),t(M'))$ for all right $\mathcal D$-modules $M, M'$) and therefore takes projectives to projectives and injectives to injectives. Therefore, for all right $\mathcal D$-modules $M, M'$ there exist functorial $k$-linear bijections 
\begin{equation}\label{bijections}
{\rm Ext}^i_{{\rm mod-}D}(M,M')\cong{\rm Ext}^i_{\mathcal D{\rm -mod}}(t(M),t(M')).
\end{equation}

In view of all this, in the case that char$k=0$ and $R=k[[x_1,\dots,x_n]]$ we are going to view all right $\mathcal D$-modules $M$ as left $\mathcal D$-modules $t(M)$ and by abuse of notation, will denote $t(M)$ by $M$. For example, if $M$ is a left $\mathcal D$-module, then $t(D(M))$ also is a left $\mathcal D$-module that we are going to denote simply $D(M)$. This convention is used in the statements of Theorems \ref{HP} and \ref{maincoro} where the left $\mathcal D$-module $t(D(R))=t(E)$ is denoted by $E$. With this notational convention and taking into account the $k$-linear bijections (\ref{bijections}) and the fact that the functors $D$ and $t$ commute, Corollary \ref{RE} in the case that $R=k[[x_1,\dots, x_n]]$ and ${\rm char}k=0$ takes the following form.

\begin{corollary}\label{REformal}
If $M$ is a left $\mathcal D$-module, there exist $k$-linear bijections $${\rm Ext}^i_{\mathcal D{\rm  -mod}}(R,D(M))\cong {\rm Ext}^i_{\mathcal D{\rm -mod}}(M,E)$$
that are functorial in $M$.
\end{corollary}

For the rest of the paper by a $\mathcal D$-module we mean a left $\mathcal D$-module. 

Let $M$ be a $\mathcal D$-module. The de Rham complex $\Omega^\bullet(M)$ of $M$ is the complex $$0\to \Omega^0(M)\to \Omega^1(M)\to \dots\to \Omega^n(M)\to 0$$ where $\Omega^s(M)=\oplus M_{i_1,\dots,i_s}$ is the direct sum of copies of $M$ indexed by all the $s$-tuples $(i_1,\dots, i_s)$ with $1\leq i_1<i_2<\dots<i_s\leq n$ and the differential $d_s:\Omega^s(M)\to \Omega^{s+1}(M)$ is given by $$d_s(x)_{i_1,\dots,i_{s+1}}=\Sigma_j(-1)^j\partial_j(x_{i_1,\dots,i_{j-1},i_{j+1},\dots,i_{s+1}})$$ for every $x\in \Omega^s(M)$ where $d_s(x)_{i_1,\dots,i_{s+1}}$ and $x_{i_1,\dots,i_{j-1},i_{j+1},\dots,i_{s+1}}$ denote the components of $d_s(x)$ and $x$ in, respectively, $M_{i_1,\dots,i_{s+1}}$ and $M_{i_1,\dots,i_{j-1},i_{j+1},\dots,i_{s+1}}$. The $i$-th de Rham cohomology group of $M$, denoted $H^i_{\rm dR}(M)$, is by definition the $i$-th cohomology group of the de Rham complex $\Omega^\bullet (M)$. The following theorem is well-known; we provide a proof for the convenience of the reader since we do not know a suitable reference.

\begin{theorem}\label{derhamext}
Let $M$ be a $\mathcal D$-module. There is a $k$-linear isomorphism $$H^i_{\rm dR}(M)\cong {\rm Ext}^i_{\mathcal D{\rm -mod}}(R, M).$$
\end{theorem}

\emph{Proof.} Let $\Delta^\bullet$ be the complex $$0\to \Delta_n\to \Delta_{n-1}\to\dots\to \Delta _1\to\Delta_0\to 0$$
where $\Delta_s=\oplus \mathcal D_{i_1,\dots,i_s}$ is the direct sum of copies of $\mathcal D$ indexed by all the $s$-tuples $(i_1,\dots, i_s)$ with $1\leq i_1<i_2<\dots<i_s\leq n$ and the differential $d_s:\Delta_s\to \Delta_{s-1}$ is given by $$d_s(x)=\Sigma_j(-1)^j(x\partial_j)_{i_1,\dots,i_{j-1},i_{j+1},\dots,i_s}$$ for every $x\in \mathcal D_{i_1,\dots,i_s}\subset\Delta_s$ where $(x\partial_j)_{i_1,\dots,i_{j-1},i_{j+1},\dots,i_s}$ denotes the element $x\partial_j\in \mathcal D_{i_1,\dots,i_{j-1},i_{j+1},\dots,i_s}\subset \Delta_{s-1}$.

The complex $\Delta^\bullet$ is a finite free resolution of $R$ in the category of left $\mathcal D$-modules. Therefore ${\rm Ext}^i_{\mathcal D{\rm -mod}}(R, M)$ is the $i$th cohomology group of the complex Hom$_{\mathcal D{\rm -mod}}(\Delta^\bullet,M)$, i.e. the complex 
$$0\to{\rm Hom}_{\mathcal D{\rm -mod}}(\Delta_0,M)\to{\rm Hom}_{\mathcal D{\rm -mod}}(\Delta_1,M)\to\dots \to{\rm Hom}_{\mathcal D{\rm -mod}}(\Delta_n,M)\to 0.$$

The map $\phi: \mathcal D\stackrel{x\mapsto x\partial_j}{\to}\mathcal D$ which is the right multiplication by $\partial_j$ is a map of left $\mathcal D$-modules. The induced map Hom$_{\mathcal D{\rm -mod}}(\phi, M)$ is nothing but the left multiplication by $\partial_j$ on $M$, i.e. $M\stackrel{x\mapsto \partial_jx}{\to}M$. Therefore the complex Hom$_{\mathcal D{\rm -mod}}(\Delta^\bullet,M)$ is nothing but the de Rham complex $\Omega^\bullet(M)$ of $M$.
\qed

\medskip

A key result of N. Switala's theory is the following \cite[5.1]{S1}.

\begin{theorem}\label{swit}
If $M$ is a holonomic $D$-module, then for all $i$ there are isomorphisms $(H^i_{\rm dR} (M))' \cong H^{n-i}_{\rm dR}(D(M))$, where $(-)'={\rm Hom}_k(-,k)$.
\end{theorem}

An important result of van den Essen \cite[Proposition 2.2]{E}\cite{S2} says

\begin{theorem}\label{vanden}
If $M$ is a holonomic $\mathcal D$-module, then $H^i_{\rm dR}(M)$ is a finite-dimensional $k$-vector space.
\end{theorem}

Corollary \ref{REformal} implies that dim$_k{\rm Ext}^i_{\mathcal D{\rm -mod}}(M,E)$ equals ${\rm dim}_k{\rm Ext}^i_{\mathcal D{\rm -mod}}(R,D(M))$ which by Theorem \ref{derhamext} equals dim$_kH^i_{\rm dR}(D(M))$ which by Theorems \ref{swit} and \ref{vanden} equals dim$_kH^{n-i}_{\rm dR}(M)$. This completes the proof of Theorem \ref{maincoro}. \qed

\end{document}